\documentclass[a4paper,12pt]{article}

\usepackage{amsfonts}
\usepackage{amssymb}

\newtheorem{theorem}{Theorem}[section]
\newtheorem{lemma}[theorem]{Lemma}
\newtheorem{corollary}[theorem]{Corollary}
\newtheorem{proposition}[theorem]{Proposition}

\newcommand{\rr}{\mathbb{R}}  %realna st
\newcommand{\cc}{\mathbb{C}}  %kompleksna st
\newcommand{\ra}{\rightarrow}  %desna puscica za preslikavo
\newcommand{\Qed}{\begin{flushright} $\Box$ \end{flushright}}
\newcommand{\proof}{\textit{Proof: }}
\newcommand{\rank}{\hbox{rank}~}
\newcommand{\im}{\hbox{Im~}}
\newcommand{\diag}{\hbox{diag~}}
\newcommand{\lin}{\hbox{lin~}}
\title{Adjacency preserving mappings on real symmetric matrices}
\author{Peter Legi\v{s}a\footnote{The author was partially supported by a grant
from the Ministry of Science of Slovenia} \\ 
\footnotesize{Department of mathematics, FMF} \\
\footnotesize{University of Ljubljana} \\
\footnotesize{Jadranska 21} \\
\footnotesize{SI-1000 Ljubljana} \\
\footnotesize{Slovenia} \\
\footnotesize{peter.legisa@fmf.uni-lj.si}
}

\begin{document}

\date{}
\maketitle

\begin{abstract}
Let $S_{n}$ denote the space of all $n\times n$ real symmetric matrices. For $n\geq 2$ we characterize maps $\Phi:S_{n}\ra S_{m}$,
which preserve adjacency, i.e. if $A,B\in S_{n}$ and $\hbox{rank}~(A-B)=1$, then $\hbox{rank}~(\Phi(A)-\Phi(B))=1$.
\end{abstract}

\vskip 24pt


\newpage

\section*{Introduction}

Wen-Ling Huang and Peter \v{S}emrl in \cite{ws} characterized 
adjacency preserving 
maps from $H_{n}$ to $H_{m}$, where $H_{n}$ 
denotes the $n\times n$ hermitian matrices over $\cc$.
They improved the results going back to Hua 
(\cite{h}, \cite{h1}). See also \cite{p1}-\cite{p9},
 \cite{w1}-\cite{w5}. 
This article considers adjacency preserving mappings 
from $S_{n}$ to $S_{m}$, where $S_{n}$ denotes the $n\times n$ symmetric matrices 
over $\rr$. The authors of \cite{ws} suggested this problem in their article. 
It turns out that the ideas and methods of their paper work in the real case 
as well
(with  modifications in some places).

The proof of the complex case uses results by Wen-Ling Huang, Roland H\"ofer
and Zhe-Xian Wan \cite{whw}, 
which hold in the real case as well. We also take advantage of 
a theorem of Alexandrov  \cite{a} on Minkowski geometries. 
Alternatively, we can use the the recent result \cite{w}
of Wen-Ling Huang  
on adjacency preserving maps from $S_{2}$ to $S_{2}$. (This paper 
 is also based on  projective geometry.)

The main result of this paper is {\bf Theorem \ref{theorem:main}}.

\section{Notation}

We will consider  only matrices over $\rr$. Let $M_{n}=M_{n}(\rr)$ be the space of all
$n\times n$ matrices. Let $S_{n}$ denote the linear subspace of all symmetric matrices
in $M_{n}$, i.e. all $A\in M_{n}$ such that $A=A^{T}$, where $A^{T}$ is the transpose 
of $A$. Let $GL(n)$ denote the group of all invertible matrices in $M_{n}$. 
Let $\hbox{lin} ~Z$ denote the real linear span of a set $Z$ (in some vector space). 
We will often look at matrices in $M_{n}$ as linear operators on  $\rr^{n}$.
So for $A\in M_{n}$, $\hbox{Im~}A=A\rr^{n}$ is the {\em image} of $A$ or the {\em column 
space} of $A$.

If we  consider $x,y\in\rr^{n}$ as $n\times 1$ matrices, $xy^{T}=x\otimes y$ 
is the rank one matrix with the property $(x\otimes y)z=\left\langle z,y\right\rangle x$ 
for $z\in\rr^n$.

If $P\in S_n$ and $P^2=P=P^T\neq 0$, then we call $P$  a {\em projection},
as it is the orthogonal projection on $\hbox{Im~}P$.
Two projections $P, Q$ are orthogonal, $P\bot Q$, iff $PQ=0$. 
If $x$ is a unit vector, then $x\otimes x$ is the projection on $\hbox{lin}~\left\{x\right\}$.

Let $e_1,\ldots e_n$ be the standard basis in $\rr^n$ and let $e_i\otimes e_j=E_{ij}$ 
be the matrix unit, i.e. the matrix with $1$ in place $(i,j)$ and zeros elsewhere.

We know that for $R,T\in M_n$, $\hbox{Im~}(R+T)\subseteq \hbox{Im~}R+\hbox{Im~}T$ 
and so \mbox{$\hbox{rank}~(R+T)$}$\leq \hbox{rank}~ R+\hbox{rank}~ T$.

For $A,B\in S_n$ let $d(A,B)=\hbox{rank}~(A-B)$. Then $(S_n,d)$ is a metric space.
We will often use 
\begin{lemma}
\label{lemma:l1}
Let $A,B,C\in M_n$ and $A+B=C$. Then $\hbox{rank}~A=\hbox{rank}~B+\hbox{rank}~C$
iff $\hbox{Im~}A=\hbox{Im~}B\oplus \hbox{Im~}C$.
\end{lemma}
Two matrices $A,B$ are {\em adjacent} if $d(A,B)=1$, i.e. $\hbox{rank}~(A-B)=1$. 
If $d(A,B)=k$, there is a sequence of consecutively adjacent matrices 
$A_0=A,A_1,\ldots , A_k=B$ (see Proposition 5.5 in \cite{wan}). Conversely,
if there is such a sequence, it is straightforward that $d(A,B)\leq k$.

Let $A,B\in S_n$ be adjacent. The line $\mathit{l}(A,B)$ joining $A$ and $B$
is the set consisting of $A,B$ and all $Y\in S_n$, which are adjacent to both $A$ and $B$.
By \cite{wan},
\[
\mathit{l}(A,B)=\left\{A+\lambda(B-A);\lambda\in\rr\right\}.
\]
If $P\in S_n$ is a projection, let $PS_nP=\left\{PAP;A\in S_n\right\}=\left\{C\in S_n;PCP=C\right\}$.

\begin{proposition}
\label{proposition:p1}
For $A,B,S\in S_n$, $R\in GL(n)$, and $c\in\rr\backslash\left\{0\right\}$ we have 
 \linebreak[4] $d(A+S,B+S)=d(A,B)=d(RAR^T,RBR^T)=d(cA,cB)$. Consequently, 
these are equivalent:
\begin{description}
\item[i)] $A$ is adjacent to $B$;
\item[ii)] $A+S$ is adjacent to $B+S$;
\item[iii)] $RAR^T$ is adjacent to $RBR^T$;
\item[iv)] $cA$ is adjacent to $cB$.
\end{description}
\end{proposition}

\begin{corollary}
\label{corollary:c1}
Let $\Phi:S_n\ra S_m$ be a map preserving adjacency, i.e. $A$ is adjacent to 
$B$ implies $\Phi(A)$ is adjacent to $\Phi(B)$. Let $\Psi(A)=\Phi(A)-\Phi(0)$
for $A \in S_n$. Then $\Psi$ is adjacency preserving and $\Psi(0)=0$.
\end{corollary}

\begin{theorem}[MAIN THEOREM]
\label{theorem:main}
Let $m$, $n$ be natural numbers, $n\geq 2$. Let $\Phi:S_n\ra S_m$ be a map preserving adjacency, with $\Phi(0)=0$. Then either:
\begin{description}
\item[i)] There is a  rank one  matrix $B\in S_m$ and a function $f:S_n\ra\rr$ such that for $A\in S_n$
\[
\Phi(A)=f(A)B.
\]
In this case we say $\Phi$ is a {\bf degenerate} adjacency preserving map.
\item[ii)] We have  $c\in\left\{-1,1\right\}$, $R\in GL(m)$ 
such that for $A\in S_n$,
\[
\Phi(A)=cR
\left[
\begin{array}{cc}
A&0  \\
0&0
\end{array}
\right]R^T.
\]
In this case we say $\Phi$ is a {\bf standard} map. (Obviously, in this case $m \geq n$.)
\end{description}
\end{theorem}

\section{Preliminary results}

We borrow Lemma 2.1. in \cite{whw}:

\begin{lemma}
\label{lemma:l2}
Let $G\in S_n$ and let $l$ be a line in $S_n$. Then either:
\begin{description}
\item[i)]There is $k$ such that $d(G,X)=k$ for all $X\in l$ or
\item[ii)]There is a point $K\in l$  such that $d(G,X)=d(G,K)+1$ for all $X\in l$, $X\neq K$.
\end{description}
\end{lemma}

\begin{lemma} 
\label{lemma:l3}
Let $A\in S_n$ be adjacent to both $R$ and $\lambda R$, where $R\in S_n$ has rank one and $\lambda\neq 1$. Then $A=\mu R$ for some $\mu\in\rr$, $\mu\neq 1,\lambda$.
\end{lemma}
\proof Since $\lambda\neq 1$, $R$ and $\lambda R$ are adjacent and $A$ is contained in the line $l(R,\lambda R)$. So $A=R+\mu'(R-\lambda R)=\mu R$ and $\mu\neq\lambda,1$.\Qed

The following lemma is slightly more general then Lemma 2.3. in \cite{ws}.

\begin{lemma}
\label{lemma:l4}
Let  $P\in M_n$ be an idempotent and $A,B\in M_n$ such that $P=A+B$ and $\hbox{rank}~P=\hbox{rank}~A+\hbox{rank}~B$. Then $A,B$ are idempotents and $AB=BA=0$.
\end{lemma}
\proof By Lemma \ref{lemma:l1}, $\im P=\im A\oplus \im B$. So if $Px=0$, 
$Ax=Bx=0$ and thus $\ker P\subset\ker A$. For $y\in\im A\subset\im P$, 
$Py=y=Ay+By$, hence $y-Ay=By$. Since $y-Ay\in \im A$, we have $By=0$. 
Thus $BA=0$ and $A^2=A$. By symmetry, $AB=0$ and $B^2=B$.\Qed

\begin{lemma}
\label{lemma:l5}
Let $P_1,P_2,\ldots,P_k \in S_n$ be mutually orthogonal rank one projections and 
$P=P_1+\ldots+P_k$. Let $\xi(1),\ldots,\xi(n)$ be an orthonormal system in $\rr^{n}$ 
such that 
 $P_i(\xi(i))=\xi(i)$ for $i=1,\ldots,k$. Then $P_i(\xi(j))=\delta_{ij}\xi(j)$. 
 Let $V$ be the orthogonal matrix defined  by $Ve_i=\xi(i)$ for $i=1,\ldots,n$, 
 so that $\xi(i)$ is the $i-$th column of $V$. Then  $V^T P_i V=E_{ii}$ for $i=1,\ldots,k$.
If $A\in PS_nP=\left\{C\in S_n|PCP=P\right\}$, then
\[
V^TAV=\left[
\begin{array}{cc}
q(A)&0  \\
0&0
\end{array}\right]
\]
where $q(A)\in S_k$. We have $q(P_i)=E_{ii}$ for $i=1,\ldots,k$ and $q(P)=E_{11}+\ldots+E_{kk}$.

The mapping $q:PS_nP\ra S_k$ is linear, bijective, and $q(AB)=q(A)q(B)$ 
for $A,B,AB\in PS_nP$.
So $q(A^2)=q(A)^2$ and $q$ is a Jordan isomorphism. It preserves the distance $d$ and 
thus adjacency. Also $q(ABA)=q(A)q(B)q(A)$ for all $A,B\in PS_nP$. All 
these properties 
are shared by the mappings $h:S_k\ra S_n$ and $q^{-1}:S_k \ra PS_nP$, where
\[
h(B)=\left[
\begin{array}{cc}
B&0  \\
0&0
\end{array}\right]
\]
and $q^{-1}(B)=Vh(B)V^T$.
\end{lemma}

\begin{lemma}
\label{lemma:l6}
Let $k,n$ be natural numbers with $3\leq k \leq n$. Let $\lambda_1,\ldots,\lambda_k$
be nonzero real numbers and $P_1,\ldots,P_k\in S_n$ mutually orthogonal rank one 
projections. Let $A=\sum_{j=1}^{k}\lambda_jP_j$.
Let $B\in S_n$ have $\rank B=\rank A=k$ and let $B$ be adjacent to 
$A - \lambda_i P_i$ for all $i$. Assume that $d(B,\lambda_i P_i)=k-1$ for all $i$. 
Then $B=A$.
\end{lemma}
\proof By Lemma \ref{lemma:l1}, 
$\im B=\im(\lambda_iP_i)\oplus\im(B-\lambda_iP_i)$. 
So $\im P_i\subset \im B$ for all $i$. 
If $P=P_1+\ldots+P_k$, then $\im P\subset \im B$ and $\rank P=k$, 
so $\im P=\im B$ and consequently $PB=B=BP$. 
Thus $A,B\in PS_n P$. 
Using notation from Lemma \ref{lemma:l5}, $q(A),q(B)\in S_k$ and $q(B)$ is adjacent to $q(A)-\lambda_iq(P_i)$, $d(q(B),\lambda_iq(P_i))=k-1$. Also, $q(P)=E_{11}+\ldots+E_{kk}=I_k$ and $q(A),q(B)$ have maximal rank as elements in $S_k$. Thus we may assume that $k=n$ and $A,B$ are invertible in $S_n$, $P_1+\ldots+P_n=I$.

Now $1=\rank(B-A+\lambda_iP_i)=\rank(A^{-1}B-I+\lambda_iA^{-1}P_i)$. But $A^{-1}=\sum \lambda_{i}^{-1}P_i$, so $\lambda_i A^{-1}P_i=P_i$. Let $C=B^{-1}A\in M_n$. Then $C\in GL(n)$ and $1=\rank(C^{-1}-(I-P_i))=\rank(I-C(I-P_i))$.

Now $I=(I-C(I-P_i))+C(I-P_i)$ and $\rank C(I-P_i)=\rank(I-P_i)=n-1$. By Lemma \ref{lemma:l3}, $C(I-P_i)$ is an idempotent.

Let $f_1,\ldots,f_n$ be an orthonormal basis of $\rr^{n}$ such that $P_i f_i=f_i$. Then for $j\neq i$, $(I-P_i)f_j=f_j$, so $Cf_j=C(I-P_i)f_j=C(I-P_i)C(I-P_i)f_j=C(I-P_i)Cf_j$. Since $C$ is invertible, 
$$(I-P_i)Cf_j=f_j \textrm{ for } j\neq i.$$

Let $Cf_j=\sum_{m=1}^{n}a_mf_m$. Then 
$(I-P_i)Cf_j=Cf_j-P_iCf_j=Cf_j-a_if_i=\sum_{m\neq i}a_mf_m=f_j$.
So $a_m=0$ for $m\neq i,j$ and $a_j=1$. Thus $Cf_j=f_j+a_if_i$. 
Since $n\geq 3$, there exists $k$, $1\leq k \leq n$, $k\neq i,j$. So $Cf_j=f_j+a_kf_k$ also. Thus $Cf_j=f_j$ for all $j$ and $C=I$. This implies $A=B$.\Qed

\begin{lemma}
\label{lemma:l7}
Let $A,B\in S_m$ and let $\rank A=1$. If $\rank(A+\lambda B)=1$ for every $\lambda\in\rr$, then $B=0$.
\end{lemma}
\proof If $\rank B\geq 2$, then there exists a nonsingular $2\times 2$ submatrix in $B$. 
For $\lambda\neq 0$, we have $\rank(A+\lambda B)=\rank( B+\frac{1}{\lambda}A)\geq 2$ for 
$\lambda$  large enough, since the chosen submatrix of $(B+\frac{1}{\lambda}A)$ 
will be nonsingular. Therefore, $\rank B\leq 1$.

If $B\neq 0$, then $B$ is adjacent to $0$. Also, $A+B$ is adjacent to $0$ and $B$, so $A+B\in l(B,0)$. Thus $A+B=\mu B$ and $A+(1-\mu)B=0$ -- a contradiction.\Qed

\begin{lemma}
\label{lemma:l8}
Let $A,B\in S_n$ have $\rank n \quad (n \geq 2)$, with $A\neq B$. 
There exists a natural number $k$ and invertible matrices $A=A_0,A_1,\ldots,A_k=B$ such that the neighbours  in this sequence are adjacent and there is a matrix $C_j\in l(A_j,A_{j+1})$ with $\rank C_j=n-1$ for $j=0,\ldots,k-1$.
\end{lemma}
\proof 
This is a consequence of Lemmas 2.5 and 2.6 in \cite{whw} and is stated in 
the proof of Lemma 3.1 in the same paper. 

\begin{lemma}
\label{lemma:l9}
Let $\Phi:S_n\ra S_m$ be an adjacency preserving map. 
Let $A,B\in S_n$ be adjacent. 
Then $\Phi(l(A,B))\subset l(\Phi(A),\Phi(B))$. The restriction of $\Phi$ to $l(A,B)$ is injective.
\end{lemma}
\proof 
If $\lambda_1\neq\lambda_2$ and $C_i=A+\lambda_i(B-A)\in l(A,B)$ ($i=1,2$), 
then $C_1$ is adjacent to $C_2$ and therefore $\Phi(C_1)$ is adjacent to $\Phi(C_2)$, 
thus $\Phi(C_1)\neq\Phi(C_2)$. \Qed

\begin{lemma}
\label{lemma:l10}
Let $\Phi:S_n\ra S_m$ ($n\geq 2$) be a map preserving adjacency and $\Phi(0)=0$.
Let $\hbox{max}\left\{\rank\Phi(A)|A\in GL(n)\right\}=k$. If $k\geq 2$ and for every singular $A\in S_n$ we have $\rank\Phi(A)<k$, then $\rank\Phi(B)=k$ for every invertible $B\in S_n$.
\end{lemma}
\proof 
Let $A,B\in S_n\cap GL(n)$ with $A\neq B$ and let $\rank\Phi(A)=k$. 
By Lemma \ref{lemma:l8}, there exists a natural number $r$ and invertible matrices 
$A=A_0,A_1,\ldots,A_r=B$ such that the neighbours in this sequence are adjacent and 
for $j=0,\ldots,r-1$ there is a matrix $C_j\in l(A_j,A_{j+1})$ with $\rank C_j=n-1$. 
Hence $\rank \Phi(C_j)<k$. Now $\rank\Phi(A)=k$, $\rank\Phi(A_1)\leq k$, $\rank\Phi(C_0)<k$.
Lemma \ref{lemma:l2} (for $G=0$) tells us that $\Phi(C_0)$ is the only point on the 
line $l(\Phi(A),\Phi(A_1))$ with rank less than $k$. Since $C_0 \neq A_1$ , 
Lemma \ref{lemma:l9} tells us that $\Phi(C_0)\neq\Phi(A_1)$. So $\rank\Phi(A_1)=k$. Proceeding in this way we find $\rank\Phi(A_j)=k$ for all $j$, so $\rank\Phi(B)=k$.\Qed

\begin{lemma}
\label{lemma:l11}
Let $\Phi:S_n\ra S_m$ ($n\geq 2$) be a map preserving adjacency. 
If there are $A,B\in S_n$ with $d(\Phi(A),\Phi(B))=n$, 
then $d(\Phi(X),\Phi(Y))=d(X,Y)$ for all $X,Y\in S_n$ and $\Phi$ is injective.
\end{lemma}
\proof 
For $n=m$ this was proved (in even greater generality) by Wen-ling Huang
(Corollary
3.1 in \cite{wh2}).

We know that $d(X,Y)=k\geq 1$ implies the existence of a sequence 
$X=X_0,X_1,\ldots,X_k=Y$ of
consecutively adjacent matrices. If $\Psi:S_n\ra S_m$ is adjacency preserving, 
the neighbours in the sequence $\Psi(X_0),\Psi(X_1),\ldots,\Psi(X_k)$ are also adjacent
and therefore $d(\Psi(X),\Psi(Y))\leq k$. So 
\[
d(\Psi(X),\Psi(Y))\leq d(X,Y)
\]
for any adjacency preserving map $\Psi$.

Now the map $\Psi$, defined by $\Psi(X)=\Phi(X+A)-\Phi(A)$ for $X\in S_n$ is adjacency 
preserving by Proposition \ref{proposition:p1} and $\Psi(0)=0$. We note 
that $\rank(\Psi(B-A))=d(\Phi(B),\Phi(A))=n$.

If $Z\in S_n$ is singular, 
$\rank(\Psi(Z))=d(\Psi(Z),\Psi(0))\leq d(Z,0)=\rank Z\leq n-1$. 
Lemma \ref{lemma:l10} tells us that $\rank(\Psi(X))=n$ for every 
$X\in S_n\cap GL(n)$. In particular, if $d(C,A)=n$, i.e. $\rank(C-A)=n$,
then $n=\rank(\Psi(C-A))=\rank(\Phi(C)-\Phi(A))=d(\Phi(C),\Phi(A))$.

Let $X,Y\in S_n$ be such that $d(X,Y)=n$. For $\lambda$ large enough, 
$d(\lambda I,A)=\rank(\lambda I-A)=n$ and $d(\lambda I,X)=n$. If we set $C=\lambda I$ 
above, we see $d(\Phi(\lambda I),\Phi(A))=n$. We may substitute $\lambda I$ for $A$, $A$ for $B$
in the previous argument and get $d(\Phi(\lambda I),\Phi(X))=n$. Repeating this procedure we get $d(\Phi(X),\Phi(Y))=n$.

We have proven that $d(X,Y)=n$ implies $d(\Phi(X),\Phi(Y))=n$. 
Suppose now $d(Z,W)=\rank(Z-W)=k<n$, with $k\geq 1$. There is $U$ orthogonal such 
that $Z-W=U(\lambda_1E_{11}+\ldots+\lambda_kE_{kk})U^T$, with $\lambda_1,\ldots,\lambda_k$ 
nonzero. Let $G=W-U(E_{k+1,k+1}+\ldots+E_{nn})U^T$. 
Then $d(G,W)=\rank(G-W)=n-k$ and $(Z-W)+(W-G)=Z-G$ is invertible. Since $\Phi$ does not 
increase the metric $d$, 
$n=d(Z,G)=d(Z,W)+d(W,G)\geq d(\Phi(Z),\Phi(W))+d(\Phi(W),\Phi(G))\geq d(\Phi(Z),\Phi(G))=n$.
So $d(Z,W)=d(\Phi(Z),\Phi(W))$.

If $\Phi(X)=\Phi(Y)$ and $X\neq Y$, then $d(X,Y)\geq 1$, so $d(\Phi(X),\Phi(Y))\geq 1$ -- a contradiction.\Qed

\begin{lemma}
\label{lemma:l12}
Let $m>n\geq 2$ and let $A_1,B_1\in S_n$ with $A_1\neq B_1$. If $A,B\in S_m$ are such that 
\[
A=\left[
\begin{array}{cc}
A_1&0  \\
0&0
\end{array}
\right],\;
B=\left[
\begin{array}{cc}
B_1&0  \\
0&0
\end{array}
\right]
\]
and $C$ is adjacent to both $A$ and $B$, then there is $C_1\in S_n$ such that
\[
C=\left[
\begin{array}{cc}
C_1&0  \\
0&0
\end{array}
\right].
\]
\end{lemma}
\proof 
The matrices $A-C$ and $C-B$ have rank one. So $A-B=(A-C)+(C-B)$ has rank 
one or two. If $A$ is adjacent to $B$, then $C$ lies one the line 
$l(A,B)$, so $C=A+\lambda(B-A)$ has the desired form.

If $A-B$ has rank two, then $\im(A-B)=\im(A-C)\oplus\im(C-B)$ by Lemma \ref{lemma:l1}.
So $\im(A-C)\subset\im(A-B)$ and $C=A-(A-C)$ has the desired form.\Qed

\begin{lemma}
\label{lemma:l13}
Let $m>n\geq 2$ and let $\Phi:S_n\ra S_m$ be an adjacency preserving map with $\Phi(0)=0$. Let
\[
\Phi(I)=\left[
\begin{array}{cc}
K&0  \\
0&0
\end{array}
\right]
\]
where $K\in S_n$ has rank $n$. Then for all $A\in S_n$,
\[
\Phi(A)=\left[
\begin{array}{cc}
A_1&0  \\
0&0
\end{array}
\right]
\]
where $A_1\in S_n$.
\end{lemma}
\proof Since $n=d(\Phi(I),\Phi(0))$, Lemma \ref{lemma:l11} tells us that $d$ preserves
the distance. Suppose $P\in S_n$ is a projection of rank one. Then $d(0,P)=1$, $d(I,P)=n-1$, so $d(\Phi(I),\Phi(P))=n-1$ and $d(0,\Phi(P))=1$. Thus
\[
n=\rank\Phi(I)=\rank\Phi(P)+\rank(\Phi(I)-\Phi(P)).
\]
By Lemma \ref{lemma:l1}, $\im\Phi(I)=\im\Phi(P)\oplus\im(\Phi(I)-\Phi(P))$, so $\im\Phi(P)\subset\im\Phi(I)$ and $\Phi(P)$ has the desired form.

If $A=\lambda P$, then $A$ lies on the line $l(0,P)$, so $\Phi(A)$ lies on the line $l(0,\Phi(P))$, so $\Phi(A)=\mu \Phi(P)$ has the desired form.

Now we use the induction on the rank of $A$. Suppose we have proved the lemma 
for all matrices of $\rank k\geq 1$. Let $\rank A=k+1$. There is $U$ orthogonal 
and nonzero numbers $\lambda_1,\ldots,\lambda_{k+1}$ such that 
$A=U(\lambda_1E_{11}+\ldots+\lambda_{k+1}E_{k+1,k+1})U^T$. 
The matrix $A$ is adjacent to $B=U(\lambda_2E_{22}+\ldots+\lambda_{k+1}E_{k+1,k+1})U^T$
and to $C=U(\lambda_1E_{11}+\ldots+\lambda_{k}E_{kk})U^T$. So $\Phi(A)$ is adjacent
to
\[
\begin{array}{ccc}
\Phi(B)=\left[
\begin{array}{cc}
B_1&0 \\
0&0
\end{array}
\right]
&
and
&
\Phi(C)=\left[
\begin{array}{cc}
C_1&0 \\
0&0
\end{array}
\right]
\end{array}
\]
where $B_1,C_1\in S_n$ and $B_1 \neq C_1$. By Lemma \ref{lemma:l11}, $\Phi(B)\neq\Phi(C)$. We use Lemma \ref{lemma:l12}.\Qed

\section{Adjacent matrices in $S_2$}

Wen-Ling Huang  proved the following result  (Corollary 2 in \cite{w}):

Let $\Phi:S_2\ra S_2$ be an adjacency preserving map. Suppose there are 
$A,  B \in S_2$ such that $\Phi(A)$ and $\Phi(B)$ are not adjacent. Then 
there are $c\in\left\{-1,1\right\}$ and $T\in GL(2)$, $S\in S_2$ 
such that $\Phi(X)=cTXT^T+S$ for $X\in S_2$. 

This Corollary implies the main result of this section, Proposition 3.5. But 
we can also proceed in a way analoguous to that in \cite{ws}.

\begin{lemma}
\label{lemma:l31}
Let $\Phi:S_2\ra S_2$ be a map such that $A$ is adjacent to $B$ iff $\Phi(A)$ is adjacent to $\Phi(B)$. Then $\Phi$ is injective.
\end{lemma}
\proof If there are $A,B\in S_2$ such that $d(\Phi(A),\Phi(B))=2$, then, by Lemma \ref{lemma:l11}, $\Phi$ is injective.

Suppose now that $d(\Phi(X),\Phi(Y))\leq 1$ for all $X,Y\in S_2$. We will show this 
is impossible. Since $E_{11}$ and $E_{22}$ are not adjacent, 
$\Phi(E_{11})$ and $\Phi(E_{22})$ are not adjacent. Therefore
$\Phi(E_{11})=\Phi(E_{22})$. Similarly, $\Phi(2E_{11})=\Phi(E_{22})$. On the other hand, $E_{11}$ is adjacent to $2E_{11}$, so $\Phi(E_{11})$ is adjacent to $\Phi(2E_{11})=\Phi(E_{11})$ -- a contradiction.\Qed

We denote by $Q$ the quadratic form on $\rr^n$, defined by 
$Q(x)=x_n^2-x_1^2-x_2^2-\ldots-x_{n-1}^2$. Then $Q(x-y)$ is the 
{\bf Lorentz separation} of $x$ and $y$. A bijective linear transformation
$L:\rr^n\ra\rr^n$ is a {\bf Lorentz transformation}
if $Q(Lx)=Q(x)$ for all $x\in\rr^n$. All Lorentz transformations on $\rr^n$ form the
{\bf Lorentz group}. A mapping $f:\rr^n\ra\rr^n$ is a 
{\bf Weyl transformation} if there are: $\alpha\in\rr\backslash\left\{0\right\}$,  a Lorentz transformation $L$ and $b\in\rr^n$ such that $f(x)=\alpha Lx+b$ for all $x\in\rr^n$.

The following theorem is due to Alexandrov \cite{a}. We quote it from Lester 
\cite{lester} p. 929, who rediscovered it.

\begin{theorem}
\label{theorem:t32}
Let $D$ be an open connected subset of $\rr^n$ and let 
$f:D\ra\rr^n$  be 
an injective mapping such that $Q(x-y)=0$ iff $Q(f(x)-f(y))=0$. Then $f$ is the restriction of conformal mapping.
\end{theorem}
Any conformal mapping on $\rr^n$ is a Weyl transformation (see \cite{lester}, p. 
929 or \cite{pop}, pp. 132-133) and that is all we will need:

\begin{corollary}
\label{corollary:c33}
Let $f:\rr^n\ra\rr^n$ be an injective mapping such that $Q(x-y)=0$
iff $Q(f(x)-f(y))=0$. Then $f$ is a Weyl transformation.
\end{corollary}

We have the linear bijection $T:\rr^3\ra S_2$, defined by 
\[
Tx=\left[
\begin{array}{cc}
x_3+x_1&x_2 \\
x_2&x_3-x_1
\end{array}
\right].
\]
Now $\det(Tx-Ty)=\det (T(x-y))=Q(x-y)$. Therefore:
\begin{equation}
\label{eqn:e31}
Tx \mbox{ is adjacent to } Ty \mbox{ iff }  x\neq y  \mbox{ and } Q(x-y)=0.
\end{equation}

The following is taken from the book \cite{rr} on Hyperbolic Geometry 
by Ramsey and Richtmyer, pp. 246-250. If $L=\left[l_{ij}\right]\in M_3$ 
is a Lorentz matrix, then $\left|\det L\right|=1$ and $\left|l_{33}\right|\geq 1$.
If $\det L=1$  and $l_{33}>1$, then $L$ is a {\bf restricted Lorentz matrix}. 
If $L$ is a restricted Lorentz matrix, then there is a matrix $P_1\in M_2$ with $\det P_1=1$ such that
\[
T(Lx)=P_1(Tx)P_1^T
\]
for all $x\in\rr^3$. Now $K=-E_{11}+E_{22}+E_{33}= K^{-1}$ is a Lorentz matrix with $\det K=-1$. For $Q=E_{12}+E_{21}\in S_2$ we have $T(Kx)=Q(Tx)Q^T$.

If $L\in M_3$ is any Lorentz matrix, then there is $r\in\left\{-1,1\right\}$ such 
that $rL$ or $LK$ or $rLK$ is a restricted Lorentz matrix. It follows that for any Lorentz matrix $L\in M_3$ we have
\begin{equation}
\label{eqn:e32}
T(Lx)=c_1P(Tx)P^T
\end{equation}
where $c_1\in\left\{-1,1\right\}$, $\left|\det P\right|=1$ and $x\in\rr^3$.

\begin{corollary}
\label{corollary:c34}
Let $\Phi:S_2\ra S_2$ be a map such that $A$ is adjacent to $B$ iff $\Phi(A)$ is adjacent to $\Phi(B)$. Then there exist $c\in\left\{-1,1\right\}$, $R\in GL(2)$ and $S\in S_2$ such that
\[
\Phi(A)=cRAR^T+S \quad (A\in S_2).
\]
\end{corollary}
\proof
We consider the mapping $f:\rr^3\ra\rr^3$, defined by
\[
f(x)=T^{-1}\Phi(Tx).
\]
By Lemma \ref{lemma:l31}, $f$ is injective. If $x\neq y$ and $Q(x-y)=0$, then $Tx$ is adjacent to $Ty$, so $\Phi(Tx)$ is adjacent to $\Phi(Ty)$, so $Q(f(x)-f(y))=0$. If $f(x)=f(y)$, then $x=y$.

If $f(x)\neq f(y)$ and $Q(f(x)-f(y))=0$, then $\Phi(Tx)$ is adjacent to $\Phi(Ty)$ 
by (\ref{eqn:e31}), so $Tx$ is adjacent to $Ty$ and $Q(x-y)=0$.

We see that $Q(x-y)=0$ iff $Q(f(x)-f(y))=0$. By Corollary \ref{corollary:c33}, 
there exist $\alpha\in\rr\backslash\left\{0\right\}$, $b\in\rr^3$ and  a Lorentz 
matrix $L\in GL(3)$ such that $f(x)=\alpha Lx+b$ for all $x\in\rr^3$, hence
\[
\Phi(Tx)=\alpha T(Lx)+Tb.
\]
By (\ref{eqn:e32}), there are $c_1\in\left\{-1,1\right\}$ and $P\in GL(n)$ such that 
\[
\Phi(Tx)=\alpha c_1P(Tx)P^T+Tb,
\]
i.e.
\[
\Phi(A)=cRAR^T+S
\]
for $A\in S_2$, where $c\in\left\{-1,1\right\}$, $R\in GL(2)$ and $S\in S_2$.\Qed

\begin{proposition}
\label{proposition:p35}
Let $\Phi:S_2\ra S_2$ be an adjacency preserving mapping. Suppose $d(\Phi(G),\Phi(H))=2$ for some $G,H\in S_2$. Then there are $c\in\left\{-1,1\right\}$, $R\in GL(2)$ and $S\in S_2$ such that
\[
\Phi(A)=cRAR^T+S.
\]
\end{proposition}
\proof By Lemma \ref{lemma:l11}, $d(\Phi(X),\Phi(Y))=d(X,Y)$ for all $X,Y\in S_2$.
So $\Phi(X)$ is adjacent to $\Phi(Y)$ iff $X$ is adjacent to $Y$. 
We use Corollary \ref{corollary:c34}.\Qed

\section{Proof of theorem 1.4}

Let $n\geq 2$ and let $\Phi:S_n\ra S_m$ be a mapping preserving adjacency, $\Phi(0)=0$. Theorem \ref{theorem:main} states that $\Phi$ is either degenerate or a standard map.

\begin{lemma}
\label{lemma:l41}
Theorem \ref{theorem:main} is true if $n=2$.
\end{lemma}
\proof If $m=1$, $\Phi$ is a degenerate map. Let $m\geq 2$. We consider two cases.

{\bf Case 1:} Let $d(\Phi(A),\Phi(B))\leq 1$ for all $A,B$. 

Then $\rank\Phi(A)\leq 1$
for all $A$. Since $E_{11}$ is adjacent to $0$, $\Phi(E_{11})$ is adjacent to 
$\Phi(0)=0$, so $\rank\Phi(E_{11})=1$. Let  $A\in S_2$. Then 
$d(\Phi(A),\Phi(E_{11}))\leq 1$. So $\Phi(A)=\Phi(E_{11})$ or $\Phi(A)$ is adjacent to
$\Phi(E_{11})$. In the latter case, if $\Phi(A)\neq 0$, then 
$\Phi(A)$ is adjacent to $0$, so $\Phi(A)\in l(0,\Phi(E_{11}))$, 
thus $\Phi(A)=\lambda \Phi(E_{11})$. So $\Phi(A)=\lambda \Phi(E_{11})$ 
in any case. Thus $\Phi$ is a degenerate map.

{\bf Case 2:} We have $A,B\in S_2$ such that $d(\Phi(A),\Phi(B))=2$.

If $m=2$, then Proposition \ref{proposition:p35} ends the proof. Let $m>2$. 
By Lemma \ref{lemma:l11}, $\Phi$ preserves the distance and is injective. So $d(\Phi(I),0)=2=\rank\Phi(I)$. Since $\Phi(I)\in S_m$, there is $U\in M_m$ orthogonal such that
\[
\begin{array}{ccc}
U\Phi(I)U^T=\left[
\begin{array}{cc}
D&0 \\
0&0
\end{array}
\right]
&
\textrm{and}
&
D=\left[
\begin{array}{cc}
\lambda_1&0 \\
0&\lambda_2
\end{array}
\right]
\end{array}.
\]
Let $\Psi(A)=U\Phi(A)U^T$ for $A \in S_2$. 
Then $\Psi$ is distance preserving and $\Psi(0)=0$. By Lemma \ref{lemma:l13},
\[
\Psi(A)=\left[
\begin{array}{cc}
\Psi_1(A)&0 \\
0&0
\end{array}
\right]
\]
where $\Psi_1(A)\in S_2$ and $\Psi_{1}(0)=0$.

Obviously, $d(\Psi(A),\Psi(B))=d(\Psi_1(A),\Psi_1(B))$. So $\Psi_1:S_2\ra S_2$ 
is distance preserving. By Proposition \ref{proposition:p35}, there are 
$c\in\left\{-1,1\right\}$ and $R\in GL(2)$ such that $\Psi_1(A)=cRAR^T$. Let
\[
W=\left[
\begin{array}{cc}
R&0 \\
0&I
\end{array}
\right]\in GL(m).
\]
Then
\[
\Psi(A)=cW\left[
\begin{array}{cc}
A&0 \\
0&0
\end{array}
\right]W^T
\]
and
\[
\Phi(A)=cU^TW\left[
\begin{array}{cc}
A&0 \\
0&0
\end{array}
\right](U^TW)^T.
\]
\Qed

\begin{lemma}
\label{lemma:l42}
Let $n\geq 2$ and let $\Phi: S_n\ra S_m$ be a map preserving adjacency, with $\Phi(0)=0$. Let
\[
\Phi(I)=\left[
\begin{array}{cc}
I_n&0 \\
0&0
\end{array}
\right]\in S_m
\]
where $I_n\in M_n$ is the identity matrix. Then we can find $U\in M_n$ 
orthogonal such that for all $A\in S_n$ we have  
\[
\Phi(A)=\left[
\begin{array}{cc}
UAU^T&0 \\
0&0
\end{array}
\right].
\]
\end{lemma}
\proof
Obviously $m\geq n$. If $m>n$, then by Lemma \ref{lemma:l13}, for all 
$A\in S_n$ we have
\[
\Phi(A)=\left[
\begin{array}{cc}
\Phi_1(A)&0 \\
0&0
\end{array}
\right],
\]
where $\Phi_1(A): S_n \to  S_n$ and $\Phi_1(I)=I$. Also $\Phi_1(0)=0$ and $\Phi_1$ 
preserves adjacency. Thus it suffices to prove the theorem for
$m=n$. We wil use induction on $n$. We know our Lemma is true 
for $n=2$ using Proposition \ref{proposition:p35}. Suppose it is valid for $n-1$,
where $n\geq 3$.

Let $P\in S_n$ be a projection with $\rank P=k$. Since $d(\Phi(I),\Phi(0))=n$, 
Lemma \ref{lemma:l11} says $d(\Phi(A),\Phi(B))=d(A,B)$ for all $A,B\in S_n$.
So $\rank\Phi(A)=\rank A$ for all $A$ and $\rank\Phi(P)=k$. Since 
$d(I,P)=n-k$, 
we have $d(I,\Phi(P))=n-k$, so  $I=\Phi(P)+R_1$, where 
$\rank R_1=n-k$. By Lemma \ref{lemma:l4}, $\Phi(P)=Q$ is a projection. So
$\Phi$ maps projections into projections of the same rank. 

Suppose $k=n-1$. By Lemma 
\ref{lemma:l5}, there is a Jordan isomorphism $q:PS_nP\ra S_{n-1}$ 
which  preserves the distance $d$. 

Since $Q$ is similar to $E_{11}+\ldots+E_{n-1,n-1}$, there is $W 
\in M_n$ orthogonal such that
\[
WQW^T=\left[
\begin{array}{cc}
I_{n-1}&0 \\
0&0
\end{array}
\right]\in M_n.
\]
We define $f_1:S_{n-1}\ra S_n$ by $f_1(B)=W\Phi(q^{-1}(B))W^T$. Then
\[
f_1(I_{n-1})=\left[
\begin{array}{cc}
I_{n-1}&0 \\
0&0
\end{array}
\right]\in S_n.
\]
We use the induction hypothesis. There is $U_1\in M_{n-1}$ orthogonal ž
such that for $B\in S_{n-1}$
\[
f_1(B)=\left[
\begin{array}{cc}
U_1BU_1^T&0 \\
0&0
\end{array}
\right]\in S_n.
\]
For $A \in PS_nP$ we have

\[
\Phi(A)=W^T\left[
\begin{array}{cc}
U_1 q(A){U_1}^T&0 \\
0&0
\end{array}
\right]W \in M_n.
\]

The mapping $\Phi$ restricted to $PS_nP$ is a Jordan isomorphism 
(in particular linear)
and if $AB=0$, then $\Phi(A)\Phi(B)=0$. 
Thus $\Phi$ maps projections in $PS_nP$ into projections 
of the same rank and 
preserves the orthogonality of projections in $PS_nP$.

Since $n\geq 3$, for any rank one projections $P_1,P_2\in S_n$ with $P_1P_2=0$ there is a projection $P$ of rank $n-1$ such that $P_1,P_2\in PS_nP$. So $\Phi_1(P_1),\Phi_1(P_2)$ are rank one projections and $\Phi(P_1)\Phi(P_2)=0$.

Thus $\Phi(E_{11}),\ldots,\Phi(E_{nn})$ are mutually orthogonal rank one 
projections. By Lemma \ref{lemma:l5}, there is $V\in M_n$ orthogonal such that 
$V\Phi(E_{ii})V^T=E_{ii}$ ($i=1,\ldots,n$). By exchanging $\Phi$ with the
map $A \mapsto V\Phi(A)V^T$ we may assume 
$\Phi(E_{ii})=E_{ii}$ for $i=1,\ldots,n$.

Let $j\neq i$ and $R=E_{ii}+E_{jj}$. Since $n\geq 3$, 
there is a projection $P$ of rank $n-1$ 
 such that $RS_nR \subset PS_nP$. By the preceding paragraph
 $\Phi(R)=\Phi(E_{ii})+ \Phi(E_{jj})=E_{ii}+ E_{jj}=R$ and for 
 $A\in RS_nR$ we have $\Phi(A)=\Phi(R)\Phi(A)\Phi(R)=R\Phi(A)R$,
 so $\Phi(A)\in RS_nR$. Also, $\Phi$ restricted to $RS_nR$ 
is linear, injective, preserving the products (if the products are in $RS_nR$).

By Lemma \ref{lemma:l5}, we have the Jordan isomorphism
$q:RS_nR\ra S_2$, such that $q(E_{ii})=E_{11}$ and $q(E_{jj})=E_{22}$.
The map $K=q\Phi q^{-1}:S_2\ra S_2$ is adjacency preserving,
$K(0)=0$, $K(I_2)=I_2$. By the induction hypothesis, 
there is $U_2\in M_2$ orthogonal such that $K(B)=U_2BU_2^T$ for $B\in S_2$.
Since $K(E_{11})=E_{11}$, $K(E_{22})=E_{22}$,  $U_2= \diag (\lambda_1,\lambda_2)$,
with $\lambda_1, \lambda_2\in\left\{-1,1\right\}$. It follows there is 
$w_{ij}\in\left\{-1,1\right\}$ such that $K(E_{12}+E_{21})=w_{ij}(E_{12}+E_{21})$, 
hence   $\Phi(E_{ij}+E_{ji})
=w_{ij}(E_{ij}+E_{ji})$ and 
$\Phi(\alpha E_{ii}+\beta(E_{ij}+E_{ji})+\gamma E_{jj})=
\alpha E_{ii}+w_{ij}\beta (E_{ij}+E_{ji})+\gamma E_{jj}$ 
for $\alpha,\beta,\gamma\in\rr$. Let $w_{ii}=1$ for all $i$.

For $b\in\rr$ let $B=E_{ii}+b(E_{ij}+E_{ji})+b^2E_{jj}\in S_n$. Then $B$ 
has rank one
and $B \in RS_nR$. If $A=\left[a_{ij}\right]\in S_{n}$ has rank one, then 
there exists
a projection $Q$ with $\rank Q =n-1 $ such that $A,B\in QS_nQ$. Since 
$\Phi$ restricted to $QS_nQ$ is a Jordan map,
$\Phi(BAB)=\Phi(B)\Phi(A)\Phi(B)$. 
Also $B=BR=RB$ and so $BAB=B(RAR)B$. But 
$RAR=a_{ii}E_{ii}+ a_{ij}(E_{ij}+E_{ji})+a_{jj}E_{jj}$.
We know that $\Phi(B)=E_{ii}+w_{ij}b(E_{ij}+E_{ji})+b^2E_{jj}$, 
so $R\Phi(B)=\Phi(B)R$ and $\Phi(B(RAR)B)=\Phi(BAB)=\Phi(B)\Phi(A)\Phi(B)=
(\Phi(B)R)\Phi(A)(R\Phi(B))$. So
\begin{equation}
\label{eq:421}
\Phi(B(RAR)B)=\Phi(B)(R\Phi(A)R)\Phi(B) \quad (b\in\rr).
\end{equation}
If $\Phi(A)=\left[a_{ij}'\right]$, 
$R\Phi(A)R=a_{ii}'E_{ii}+a_{ij}'(E_{ij}+E_{ji})+a_{jj}'E_{jj}$. Equation 
\ref{eq:421} implies $a_{ii}'=a_{ii}$, $a_{jj}'=a_{jj}$ and  
$a_{ij}'=w_{ij}a_{ij}$. So for all $i,j$
\begin{equation}
\label{eq:422}
a_{ij}'=w_{ij}a_{ij}.
\end{equation}

Suppose now $T \in S_n$ is such that $t_{ij}=1$ for all $i,j$. Then $T$ 
has rank one and consequently 
$\Phi(T)=\left[w_{ij}\right]\in S_n$ has rank one. 
There exists $\lambda\in\rr$ such that $\Phi(T)=\lambda Q_2$, 
where $Q_2$ is a rank one projection. There exists a unit vector 
$x\in\rr^n$ such that $Q_2=x\otimes x$. 
Now $1=w_{11}=\left\langle \Phi(T)e_1,e_1\right\rangle=
\lambda\left\langle Q_2 e_1,e_1\right\rangle=\lambda \left\|Q_2 e_1\right\|^2$. 
So $\lambda>0$. Therefore, if $y=x\sqrt{\lambda}$, 
then $\Phi(T)=y\otimes y=y^Ty$, so $w_{ij}=y_iy_j$ 
for all $i,j$. But $w_{ii}=y_i^2=1$, so $y_i\in\left\{-1,1\right\}$ 
for all $i$. Therefore, if $V_2=\diag (y_1,\ldots,y_n)$, $V_2^T=V_2$ 
is orthogonal: $V_2^2=I$ and $\Phi(A)=V_2AV_2$. 
Thus $V_2\Phi(A)V_2=A$ for all $A\in S_n$ with $\rank A=1$. By exchanging
$\Phi$ with the map $A \to V_2\Phi(A)V_2^T$ we may assume $\Phi(A)=A$
for all $A\in S_n$ with $\rank A=1$.

Suppose $B\in S_n$ has rank less than $n$. Then $B=\sum_{i=1}^{n-1}\lambda_iP_i$, 
where $\lambda_i\in\rr$ and $P_i$ are mutually orthogonal rank one projections. 
Let $P=\sum_{i=1}^{n-1}P_i$. Since $\Phi$, restricted to $PS_nP$ is linear,

\[
\Phi(B)=\sum_{i=1}^{n-1}\Phi(\lambda_iP_i)=\sum_{i=1}^{n-1}\lambda_iP_i=B.
\]
Let $C\in S_n$ be invertible. Again, $C=\sum_{i=1}^{n}\alpha_iQ_i$, 
where $\alpha_i\neq 0$ and $Q_i$ are mutually orthogonal rank one projections.
Let $G=C-\alpha_iQ_i$. Then $\rank G=n-1$. Also $G$ is adjacent to $C$, 
so the same is true for $\Phi(C)$ and $\Phi(G)=G=C-\alpha_iQ_i$.
Since $\Phi$ preserves the distance, we have 
$n-1=\rank G=d(C,\alpha_iP_i)=d(\Phi(C),\alpha_iQ_i)$.
Since $\rank\Phi(C)=\rank C=n$, Lemma \ref{lemma:l6} implies $\Phi(C)=C$.

\Qed

\begin{lemma}
\label{lemma:l43}
Let $\Phi:S_n\ra S_m$ ($m,n\geq 3$) be an adjacency preserving map and $\Phi(0)=0$.
Suppose that for every projection $P\in S_n$ with $\rank P = n-1$ there 
is a rank one projection $Q$ such that $\Phi(PS_nP)\subset\rr Q$. Then $\Phi$ is a 
degenerate adjacency preserving map.
\end{lemma}
\proof Let $P, P_1\in S_n$ be  projections of  rank $ n-1$.  
Let $\Phi(PS_nP)\subset\rr Q$ and $\Phi(P_1S_nP_1)\subset\rr Q_1$, 
where $Q, Q_1$ are rank one projections. There is a projection $R$ 
of rank one such that $R\in PS_nP\cap P_1S_nP_1$. Since $R$ is adjacent 
to $0$, $\Phi(R)$ is adjacent to $0$, so $\Phi(R)=\lambda Q=\mu Q_1$ 
with $\lambda,\mu\neq 0$. Thus $Q=Q_1$ and $\Phi(B)\in \rr Q$ for all
$B\in S_n$ with $\rank B\leq n-1$. There is an orthogonal matrix $V$ 
such that $VQV^T=E_{11}$. Exchanging $\Phi$ for the map
$X\longmapsto V\Phi(X)V^T$ we may assume $Q=E_{11}$. So $\Phi(B)
\in \rr E_{11}$ for all $B$ with $\rank B\leq n-1$.

If $A\in S_n$ is invertible, then 
\begin{equation}
\label{eqn:411}
A=\sum_{j=1}^{n}\lambda_jP_j,
\end{equation}
where $\lambda_j$ are nonzero and $P_j$ are mutually orthogonal 
rank one projections. So $A$ is adjacent to $B=\sum_{j=2}^{n}\lambda_jP_j$.
Thus $\Phi(A)$ is adjacent to $\Phi(B)=\lambda E_{11}$ and 
$\rank\Phi(A)\leq 2$.

{\bf Case 1:} Assume $\rank \Phi(A) \leq 1$ for all 
$A \in S_n \cap GL(n)$. 

We claim $\Phi(S_n \cap GL(n)) \subset \rr E_{11}$.
Suppose, on the contrary, that there exists $A$ invertible such that 
$\Phi(A)=Z \notin \rr E_{11}$. Then $\rank Z =1$. Let
$A=\sum_{j=1}^{n}\lambda_jP_j$ as in (\ref{eqn:411}). Let 
$B=\sum_{j=2}^{n}\lambda_jP_j$. Then $Z$ is adjacent to 
$\Phi (B)=\lambda E_{11}$ and to 0.
If $\lambda \neq 0$, then $Z$ lies on the line $l(0,\lambda E_{11})$ --
a contradiction. So $\Phi(B)=0$. By Lemma \ref{lemma:l9}, $\Phi$
maps the line $l(B,A)=\left\{B+\lambda P_1;\lambda \in \rr  \right\}$
into the line $l(0,Z)= \rr Z$ injectively. So there is 
$\lambda \in \rr, \quad \lambda \neq  \lambda _1$ such that
$\Phi (B+\lambda P_1)=\frac{1}{2}Z$.
Now $C_1=\sum_{j=2}^{n-1}{\lambda}_j P_j + {\lambda}_1 P_1$ is adjacent
to $A$ and $C_2=\sum_{j=2}^{n-1}{\lambda}_j P_j + \lambda P_1$
is adjacent to $B+\lambda P_1$, so $B_1= \Phi(C_1)$
is adjacent to $Z$, $B_2= \Phi(C_2)$ is adjacent to $\frac{1}{2}Z$. Both 
$B_1$ and $B_2$ are in $\rr E_{11}$. Since $Z$ is adjacent to $\frac{1}{2}Z$, 
we have $B_1, B_2 \in l(Z, \frac{1}{2} Z) \subset \rr Z$. So $B_1, B_2=0$ - a
contradiction
with the fact that $B_1, B_2$ are adjacent.
 Thus, if $\rank\Phi(A)\leq 1$ for all invertible 
$A\in S_n$, the proof is finished.

{\bf Case 2:}  Suppose there is $A\in S_n\cap GL(n)$ with $\rank\Phi(A)=2$.

By Lemma \ref{lemma:l10},
$\rank\Phi(X)=2$ for every $X\in S_n\cap GL(n)$. So $\Phi(I)$ 
and $\Phi(I+E_{11})$ have 
rank two and are adjacent. Let $D=E_{22}+\ldots+E_{nn}$. Then $D$ is
adjacent to $I$ and to 
$I+E_{11}$, implying that $\Phi(D)$ is adjacent to $\Phi(I)$. 
Thus  $\Phi(D) \neq 0$ and
 $\Phi(D)=\lambda E_{11}$ for some $\lambda\neq 0$. Since $\Phi(D)$ is
adjacent to $\Phi(I)$ and to  $\Phi(I+E_{11})$, we have $\Phi(I)=\lambda E_{11}+K$ 
and $\Phi(I+E_{11})=\lambda E_{11}+K'$, with $\rank K=\rank K'=1$. We claim that 
\begin{equation}
\label{eqn:431}
\rank(K+\mu E_{11})=\rank(K'+\mu E_{11})=2 \quad \textrm{for} \quad  \mu\neq 0.
\end{equation}
In fact, if $\mu\neq 0$ and $\rank(K+\mu E_{11})=1$, then $K+\mu E_{11}$ is adjacent
to $0$ and to $\mu E_{11}$, so $K+\mu E_{11}$ is on the line $l(0,\mu E_{11})$, so 
$K+\mu E_{11}=\gamma\mu E_{11}$ -- a contradiction
to the fact that $\rank \Phi(I)=2$.

Now $\Phi(E_{11}+\sum_{j=3}^{n} E_{jj})=\lambda' E_{11}$ and is adjacent to 
$\Phi(I)=\lambda E_{11}+K$. Thus $\rank((\lambda-\lambda')E_{11}+K)=1$, 
which implies by (\ref{eqn:431}) that $\lambda=\lambda'$. Also $\Phi(2E_{11}+\sum_{j=3}^nE_{jj})=\lambda'' E_{11}$ is adjacent to $\Phi(I+E_{11})=\lambda E_{11}+K'$. As before, $\lambda''=\lambda$. But $E_{11}+\sum_{j=3}^nE_{jj}$ is adjacent to $2E_{11}+\sum_{j=3}^nE_{jj}$, so $\lambda E_{11}$ is adjacent to $\lambda E_{11}$ -- a contradiction.

So $\rank\Phi(A)\leq 1$ for all invertible $A\in S_n$ and the proof is finished.\Qed

\begin{lemma}
\label{lemma:l44}
Let $\Phi:S_n\ra S_m$ ($m,n\geq 3$) be an adjacency preserving map with $\Phi(0)=0$. Assume that for every projection $P$ with $\rank P=n-1$ the restriction of $\Phi$ to $PS_nP$ is a standard map. Then $\Phi$ is a standard adjacency preserving map.
\end{lemma}
\proof Let $D=E_{11}+\ldots+E_{n-1,n-1}$. There are $c\in\left\{-1,1\right\}$ and $T\in GL(m)$ such that for $B\in DS_nD$,
\[
\Phi(B)=cT\left[
\begin{array}{cc}
B&0  \\
0&0
\end{array}
\right]T^T.
\]
If $\Psi(X)=cT^{-1}\Phi(X)(T^{-1})^T$ for $X\in S_n$, then $\Psi(0)=0$, $\Psi$
preserves adjacency and
\begin{equation}
\label{eq:441}
\Psi(B)=\left[
\begin{array}{cc}
B&0  \\
0&0
\end{array}
\right] \quad (B\in DS_nD).
\end{equation}
In particular, $\Psi(E_{11})=E_{11}$. If $Q$ is a rank one projection, then
we claim 
$\Psi(Q)\geq 0$. There exists a projection $P$ of rank $n-1$, such that
$E_{11}, Q\in PS_nP$. The restriction of $\Phi$ to $PS_nP$ is a standard map. This is also true for the restriction of $\Psi$ to $PS_nP$. A standard map $\Omega$ has either the property $X\geq 0$ implies $\Omega(X)\geq 0$ or $X\geq 0$ implies $-\Omega(X)\geq 0$. Since $\Psi(E_{11})\geq 0$, $\Psi(Q)\geq 0$.

Thus $\Psi(E_{nn})\geq 0$ and $\Psi(E_{nn})$ is adjacent to $0$, so 
$\rank\Psi(E_{nn})=1$. So $\Psi(E_{nn})=s x\otimes x$ for some unit vector 
$x\in\rr^m$ and $s>0$. We show that $x\notin \lin\left\{e_1,\ldots,e_{n-1}\right\}$ 
(This implies $m\geq n$.). If $x=\alpha_1e_1+\ldots+\alpha_{n-1}e_{n-1}\in \rr^m$ and 
$y=\alpha_1e_1+\ldots+\alpha_{n-1}e_{n-1}\in\rr^n$, then by (\ref{eq:441})
we have  
$\Psi(sy\otimes y)=sx\otimes x$. There exists a projection $P_1$ of rank $n-1$ 
such that $E_{nn}, sy\otimes y\in P_1S_nP_1$. But $\Psi(E_{nn})=\Psi(sy\otimes y)$.
Since the restriction of $\Psi$ to $PS_nP$ is standard and thus injective, this is a contradiction.

We construct an invertible $R\in M_m$ such that $Re_i=e_i$ for $i=1,\ldots,n-1$ and $Rx=s^{-\frac{1}{2}}e_n$. Now $Rs(x\otimes x)R^T=sRx\otimes Rx=e_n\otimes e_n=E_{nn}$. For $i,j\leq n-1$ we have $RE_{ij}R^T=R(e_i\otimes e_j)R^T=(Re_i)\otimes(Re_j)=e_i\otimes e_j=E_{ij}$. We define $\Phi_1:S_n\ra S_m$ by $\Phi_1(X)=R\Psi(X)R^T$. Then (\ref{eq:441}) is true if we replace $\Psi$ by $\Phi_1$. Also $\Phi_1(E_{nn})=E_{nn}$. 

Let $R_i=I-E_{ii}\in S_n$. Then $R_n=D$ and $\Phi_1(R_n)=E_{11}+\ldots+E_{n-1,n-1}=E-E_{nn}$, where $E=E_{11}+\ldots+E_{nn}\in S_m$. The restriction of $\Phi_1$ to $R_iS_nR_i$ is a standard map and thus linear. So $\Phi_1(R_i)=\Phi_1(E_{11})+\ldots+\Phi_1(E_{i-1,i-1})+\Phi_1(E_{i+1,i+1})+\ldots+\Phi_1(E_{nn})=E-E_{ii}$ for $i=1,\ldots,n$.

Since $I$ is adjacent to $R_i$, $\Phi_1(I)$ is adjacent to $E-E_{ii}$ for all $i$. Thus
\begin{equation}
\label{eq:442}
E-E_{ii}=\Phi_1(I)+T_i, \mbox{ with } \rank T_i=1.
\end{equation}
Thus $\rank\Phi_1(I)\geq n-2$. If $\rank\Phi_1(I)=n-2$, then 
$\rank(E-E_{ii})=\rank\Phi_1(I)+\rank T_i$, so 
$\im(E-E_{ii})=\im\Phi_1(I)\oplus\im T_i$ by Lemma \ref{lemma:l1} 
and $\im \Phi_1(I)$ is a subspace  in $\lin(\left\{e_1,\ldots,e_n\right\}\backslash\left\{e_i\right\})$ for all $i$. Thus $\im\Phi_1(I)=\left\{0\right\}$ and $\im(E-E_{ii})=\im T_i$ -- a contradiction, since $n\geq 3$.

Suppose $\rank\Phi_1(I)=n-1=\rank(E-E_{ii})$. By (\ref{eq:442}) we have 
$\Phi_1(I)=(E-E_{ii})-T_i$ with $\rank T_i=1$. Let $T_i=\lambda_i(y_i\otimes y_i)$
with $y_i$ a unit vector. If $y_i\notin \im(E-E_{ii})$ , then $\rank\Phi_1(I)=n$ 
and that is 
a contradiction. So $y_{i}\in\im(E-E_{ii})$ and $\im\Phi_1(I)\subset\im(E-E_{ii})$ 
for all $i$. Thus once again $\Phi_1(I)=\left\{0\right\}$ -- a contradiction.

Thus $\rank\Phi_1(I)=n$.   Now $\Phi_1(I)$ is adjacent to $E-E_{ii}$ for all $i$. Also $n=d(\Phi_1(I),\Phi_1(0))=\rank\Phi_1(I)$. By Lemma \ref{lemma:l11}, $d(I,E_{ii})=n-1=d(\Phi_1(I),\Phi_1(E_{ii}))=d(\Phi_1(I),E_{ii})$. By Lemma \ref{lemma:l6}, $\Phi_{1}(I)=E$.

By Lemma \ref{lemma:l42}, we can find an orthogonal matrix $U\in M_n$ such that for $A\in S_n$
\[
\Phi_1(A)=\left[
\begin{array}{cc}
UAU^T&0  \\
0&0
\end{array}
\right].
\]
So $\Phi_1$ is a standard map and therefore $\Phi$ is a standard map.\Qed

\begin{lemma}
\label{lemma:l45}
The statement of Theorem \ref{theorem:main} is true for $n=3$.
\end{lemma}
\proof Let  $P\in S_3$ be any projection of rank 2.  By \ref{lemma:l41},
the mapping $\Phi$ restricted to $PS_3P$ is 
either standard or degenerate. 
If $\Phi$ restricted to $PS_3P$ is degenerate for all projections 
$P\in S_3$ of rank 2, 
Lemma \ref{lemma:l43} tells us that $\Phi$ is degenerate. If $\Phi$ 
restricted to 
$PS_3P$ is standard for all such $P$, then 
Lemma \ref{lemma:l44} tells us that $\Phi$ 
is a standard map.

Suppose there exist two projections $P$ and $Q$ of rank 2 such that $\Phi$ restricted
to $PS_3P$ is degenerate and $\Phi$ restricted to $QS_3Q$ is standard. Then $m\geq 2$. 
If $R\in S_3$ has rank one, then $R$ is adjacent to $0$, so $\Phi(R)$ is adjacent to 
$\Phi(0)=0$ and has rank one. There exists a rank one matrix $R_1\in QS_3Q$ such that 
the rank one matrices $\Phi(R),\Phi(R_1)$ are linearly independent. (If this is not true,
then $\Phi(R_1)=\lambda(R_1)\Phi(R)$ for all $R_1\in QS_3Q$. Since $\Phi$ restricted 
to $QS_3Q$ is standard and $\rank Q=2$, this is impossible.) There exists 
a rank two projection $R_2$ such that $R,R_1\in R_2S_3R_2$. Then since $\Phi(R)$ 
and $\Phi(R_1)$ are linearly independent, $\Phi$ restricted to $R_2S_3R_2$ is not
degenerate. Hence it is standard and therefore real linear. Thus for any rank one operator $R\in S_3$ we have
\begin{equation}
\label{eq:451}
\Phi(\lambda R)=\lambda\Phi(R) \quad (\lambda\in\rr).
\end{equation}
Let $T\in S_3$. We define $\Phi_T:S_3\ra S_m$ by $\Phi_T(X)=\Phi(X+T)-\Phi(T)$. 
Then $\Phi_T(0)=0$ and $\Phi_T$ is an adjacency preserving map by Proposition 
\ref{proposition:p1}. 

We show that $\Phi_T$ is neither standard nor degenerate. If $\Phi_T$ was standard, 
then $\Phi_T$  is real linear, so 
$\Phi(Y)=\Phi((Y-T)+T)=\Phi_T(Y-T)+\Phi(T)=\Phi_T(Y)-\Phi_T(T)+\Phi(T)$. Letting 
$Y=0$ we get $0=\Phi(0)=\Phi_T(0)+\Phi(T)-\Phi_T(T)=\Phi(T)-\Phi_T(T)$. So 
$\Phi=\Phi_T$ is standard -- a contradiction.

If there exists a rank one operator $G$ such that $\Phi_T(X)\in \rr G$ for all 
$X\in S_3$, then for $Y\in S_3$ we have $\Phi(Y)=\Phi_T(Y-T)+\Phi(T)=\Phi(T)+\lambda(Y)G$. 
Thus $0=\Phi(T)+\lambda(0)G$ and $\Phi(Y)=(\lambda(Y)-\lambda(0))G$ for all $Y\in S_3$
and $\Phi$ is degenerate -- a contradiction. 

As in the beginning of the proof of the Lemma, there are rank two projections 
$P_T$ and $Q_T$ such that the restriction of $\Phi_T$ to $P_TS_3P_T$ is degenerate and 
the restriction of $\Phi_T$ to $Q_TS_3Q_T$ is standard. If $R$ is a rank one matrix in
$S_3$, then by (\ref{eq:451}) $\Phi_T(\lambda R)=\lambda \Phi_T(R)$ for $\lambda\in\rr$. 
So $\Phi(\lambda R+ T)-\Phi(T)=\Phi_T(\lambda R)=\lambda \Phi_T(R)=
\lambda(\Phi(R+T)-\Phi(T))$, i.e.
\begin{equation}
\label{eq:452}
\Phi(\lambda R+T)=\Phi(T)+\lambda(\Phi(R+T)-\Phi(T)).
\end{equation}
Now we will prove that if $A_1,A_2,\ldots,A_p\in S_3$ have rank one, then
\[
\Phi(A_1+A_2+\ldots+A_p)=\Phi(A_1)+\ldots+\Phi(A_p)
\]
by induction on $p$. It is true for $p=1$. Assume it holds for $p$. 
Let $A_1,\ldots,A_{p+1}\in S_3$ have rank one. Then
\[
\Phi(A_1+\ldots+A_p+\lambda A_{p+1}) =    \Phi(A_1+\ldots+A_p)+\lambda(\Phi(A_1+\ldots+A_{p+1})-\Phi(A_1+\ldots+A_p))
\]
by (\ref{eq:452}), so by the induction hypotesis, 
\[
\Phi(A_1+\ldots+A_p+\lambda A_{p+1})  =  \Phi(A_1)+\ldots+\Phi(A_p)+\lambda\Phi(A_1+\ldots+A_{p+1})-\lambda(\Phi(A_1)+\ldots+\Phi(A_p)).
\]
Since $A_2+\ldots+A_p+\lambda A_{p+1}$ is adjacent to $A_1+A_2+\ldots+A_p+\lambda A_{p+1}$, we have $\Phi(A_1+\ldots+A_p+\lambda A_{p+1})$ is adjacent to $(\Phi(A_2)+\ldots+\Phi(A_p)+\lambda\Phi(A_{p+1}))$. Thus
\begin{eqnarray}\nonumber
\Phi(A_1+\ldots+A_p+\lambda A_{p+1})&-&\Phi(A_2)-\ldots-\Phi(A_p)-\lambda\Phi(A_{p+1})= \\  \nonumber 
=\Phi(A_1)+\lambda\Phi(A_1+\ldots+A_{p+1})&-&\lambda\Phi(A_1)-\ldots-\lambda\Phi(A_p)-\lambda\Phi(A_{p+1})= \\  \nonumber
=\Phi(A_1)+\lambda(\Phi(A_1+\ldots+A_{p+1})&-&(\Phi(A_1)+\ldots+\Phi(A_p)+\Phi(A_{p+1})))
\end{eqnarray}
has rank one for all $\lambda\in\rr$. By Lemma \ref{lemma:l7},
\[
\Phi(A_1+\ldots+A_p+A_{p+1})=\Phi(A_1)+\ldots+\Phi(A_{p+1}).
\]
If $A\in S_3$, then $A=\sum_{i=1}^3\lambda_iP_i$, where $P_i\in S_3$ are rank one projections. So $\Phi(A)=\sum_{i=1}^3\lambda_i\Phi(P_i)$. It follows that $\Phi$ is linear.

Now $\Phi$ maps the rank two operator $P$ into an operator of rank at most 1. By Lemma \ref{lemma:l11}, $\rank\Phi(A)\leq 2$ for all $A\in S_3$.

Let $\left\{f_1,f_2\right\}\subset\rr^3$  be an orthonormal system such
that $Qf_i=f_i$ for $i=1,2$. There exists $U\in M_3$ orthogonal such 
that $Ue_i=f_i$ for $i=1,2$. Then $QUe_i=Ue_i$ and $U^TQUe_i=e_i$ 
for $i=1,2$. Since $U^TQU\in S_3$ has rank two, we have 
$U^TQU=E_{11}+E_{22}=E_2$. If $A\in E_2S_3E_2$, then 
$U^TQUAU^TQU=A$, so $Q(UAU^T)Q=UAU^T$, so $UAU^T\in QS_3Q$. 
Now $\Phi$ restricted to $QS_3Q$ is standard.
So there are $c\in\left\{-1,1\right\}$ and $T$ invertible in $M_m$ such that
\[
\Phi(UAU^T)=cT\left[
\begin{array}{cc}
UAU^T&0  \\
0&0
\end{array}
\right]T^T.
\]

Therefore we may assume that for 
$A\in E_2S_3E_2$ we have
\[
\Phi(A) = \left[\begin{array}{cc}A&0 \\0&0\end{array}\right] =h(A) \in S_m.
\]
If $F$ is any rank two projection in $S_3$, the restriction of $\Phi$ to $FS_3F$ is 
either standard or degenerate. (Look at the beginning of the proof of this Lemma.) 
The matrix $\Phi(E_{33})$ is adjacent to $0$. Hence $\Phi(E_{33})=sx\otimes x$,
where $s\neq 0$ and $x$ is a unit vector. If $x\notin \lin\left\{e_1,e_2\right\}$,
then $\Phi(I)=E_{11}+E_{22}+sx\otimes x$ has rank 3. But $\rank\Phi(I)\leq 2$. 
So $x\in \lin\left\{e_1,e_2\right\}$ and hence $\Phi(E_{33})\in E_2S_3E_2$.

There exists a rank one projection $R_1\in E_2S_3E_2$ such that 
\[
\Phi(R_1)= h(R_1)=\left[
\begin{array}{cc}
R_1&0  \\
0&0
\end{array}
\right]
\]
and $\Phi(E_{33})$ are linearly independent.
There is a projection $R_2\in S_3$ 
of rank 2 such that $R_1,E_{33}\in R_2S_3R_2$. The restriction of $\Phi$ to 
$R_2S_3R_2$ is standard, since $\Phi(R_1)$ and $\Phi(E_{33})$ are 
linearly independent.

Since $\Phi(R_1)=h(R_1) \geq 0$, $\Phi(E_{33})\geq 0$. So
\[
\Phi(E_{33})=\left[
\begin{array}{cc}
cP_2&0  \\
0&0
\end{array}
\right]
\]
where $c>0$ and $P_2\in S_2$ is a rank one projection. Let $U_1\in M_2$ be an orthogonal matrix such that $U_1^TP_2U_1=E_{22}$. We define matrices $G\in M_3$, $V\in M_m$ by
\[
G=\left[
\begin{array}{cc}
U_1&0 \\
0&c^{-\frac{1}{2}}
\end{array}
\right],
V=\left[
\begin{array}{cc}
U_1&0 \\
0&I
\end{array}
\right].
\]
Then $GE_{33}G^T=c^{-1}E_{33}$, so $\Phi(GE_{33}G^T)=\left[\begin{array}{cc}P_2&0 \\0&0\end{array}\right]$. 

If we define $\Theta(X)=V^T\Phi(GXG^T)V$ for $X\in S_3$, then once more 
$\Theta:S_3\ra S_m$ is a linear adjacency preserving map with $\Theta(E_{33})=E_{22}$.
If $A\in E_2S_3E_2$, then $GAG^T\in E_2S_3E_2$, so 

$$\Theta(A)=h(A)= \left[\begin{array}{cc}A&0 \\0&0\end{array}\right].$$ 

If $P_1$ is a rank 
two projection in $S_3$, then, as before, $\Theta$ restricted to $P_1S_3P_1$ is either 
standard or degenerate. Now $\Theta(E_{22}+E_{33})=2E_{22}$, so $\Theta$ restricted to $(E_{22}+E_{33})S_2(E_{22}+E_{33})$ is degenerate. Therefore, $\Theta(E_{23}+E_{32})=\alpha E_{22}$, with $\alpha\neq 0$.

Since $\Theta(E_{11}+E_{33})=E_{11}+E_{22}$ has rank two, 
the restriction of $\Theta$ to $(E_{11}+E_{33})S_3(E_{11}+E_{33})$ 
is a standard map. As before, there are $c_1\in\left\{-1,1\right\}$
and $W_1\in GL(m)$ such that for $A\in (E_{11}+E_{33})S_3(E_{11}+E_{33})$ 
we have
\[
\Theta(A)=c_1W_1\left[
\begin{array}{cc}
A&0  \\
0&0
\end{array}
\right]W_1^T.
\]
But $\Theta(E_{11})=E_{11}$ and $\Theta(E_{33})=E_{22}$. So
\[
c_1W_1(e_1\otimes e_1)W_1^T=c_1(W_1e_1)\otimes(W_1e_1)=e_1\otimes e_1.
\]
This implies $c_1=1$ and $We_1=\pm e_1$. By exchanging $W$ with $-W$ if necessary we may assume $We_1=e_1$. Similarly, $We_3=de_2$, where $d\in\left\{-1,1\right\}$. This implies
\[
\Theta(E_{13}+E_{31})=W(e_1\otimes e_3+e_3\otimes e_1)W^T=We_1\otimes We_3+We_3\otimes We_1=d(E_{12}+E_{21}).
\]
Let $A=\left[1,1,1\right]^T\left[1,1,1\right]=E_{11}+E_{22}+E_{33}+(E_{12}+E_{21})+(E_{13}
+E_{31})+(E_{23}+E_{32})$. Since $A$ has rank $1$, $A$ is adjacent to $0$, so $\Theta(A)$ 
is adjacent to $0$ and has rank one. We calculate 
$\Theta(A)=E_{11}+(2+\alpha)E_{22}+(1+d)(E_{12}+E_{21})$ and 
$\det\Theta(A)=2+\alpha-(1+d)^2=\alpha-2d=0$, since $d^2=1$. So $\Theta(E_{23}+E_{32})=2dE_{22}$.

Let now $B=\left[0,d,-1\right]^T\left[0,d,-1\right]=d^2E_{22}+E_{33}-d(E_{23}+E_{32})$. 
Then $B$ has rank one and is adjacent to $0$. 
But $\Theta(B)=(1-d^2)E_{22}=0$ -- 
a contradiction.\Qed

\medskip

{\em End of proof of theorem \ref{theorem:main}}

\medskip

Let $n\geq 4$. Our induction hypothesis is that every adjacency preserving 
and zero preserving map from $S_k$ to $S_m$ ($2\leq k<n$) is either standard 
or degenerate.  Let $\Phi:S_n\ra S_m$ be an adjacency preserving map and $\Phi(0)=0$.

Let $P\in S_n$ be a projection of rank $n-1$. By Lemma \ref{lemma:l5} we 
know that $PS_nP$ is isomorphic to $S_{n-1}$. By the assumption,
$\Phi$ restricted to $PS_nP$ is either standard or degenerate. Let $Q\in S_n$ be another projection of rank $n-1$. There exists a projection $R$ of rank $n-2\geq 2$ with $PR=QR=R$, so that $R\in PS_nP\cap QS_nQ$. If $\Phi$ restricted to $PS_nP$ is degenerate, then $\Phi$ restricted to $RS_nR$ is degenerate, hence $\Phi$ restricted to $QS_nQ$ cannot be standard and is thus degenerate. By Lemma \ref{lemma:l43}, $\Phi$ is degenerate.

So if $\Phi$ restricted to $PS_nP$ is degenerate, then $\Phi$ is degenerate.
The remaining possibility is that $\Phi$ restricted to $PS_nP$ is standard.
Then obviously this is true if we replace $P$ by  $Q$.
By Lemma \ref{lemma:l43}, $\Phi$ is standard.

\section*{Acknowledgement}
The author thanks Professor Peter \v{S}emrl for suggesting this problem 
to him and for providing the preprint of the article \cite{ws}. He also 
thanks Professor Wen-ling Huang for preprints of \cite{w} and \cite{wh2}.


\begin{thebibliography}{99}
 \bibitem{ws} W.-L. Huang, P. \v{S}emrl: Adjacency preserving maps on 
 hermitian matrices, preprint (2006), to appear in Canad. J. Math.
\bibitem{h} L.K.~ Hua, Geometries of matrices $I, I_1, II, III$, 
Trans. Amer. Math. Soc. 57
(1945) 441--481, 57(1945) 482--490, 61(1947) 193--228, 61(1947) 229--255.

\bibitem{h1} L.K.~ Hua, Geometries of symmetric matrices over any field with
characteristic other than two, Ann. of Math. (2) 50(1949) 8--31.
\bibitem{whw} W.-L. Huang, R. H\"ofer, Z.-L. Wan, Adjacency preserving mappings of symmetric and hermitian matrices,
		Aequationes Math. 67 (2004) 132--139.
\bibitem{a} A.D.~ Alexandrov, On the axioms of relativity theory, 
Vestnik Leningrad. Univ. Math. 19(1976) 5--28. 
\bibitem{lester} J.A.~ Lester, Distance preserving transformations, in: 
F. Buekenhout (Eds.), Handbook of Incidence Geometry, Elsevier,
Amsterdam 1995,
921--944. 
\bibitem{w}  W.-L. Huang, Adjacency preserving mappings of $2\times 2$ 
 Hermitian matrices, Hamburger Beitr\"age zur Mathematik Heft 235, Universit\"at 
 Hamburg 2006, to appear in Aequationes Math.
 
 \bibitem{wan} Z.-X. Wan, Geometry of matrices,World Scientific, Singapore, 1996.
 
\bibitem{rr} A.Ramsey, R.D. Richtmyer, Introduction to hyperbolic geometry, Universitext 75, Springer Verlag,
		New York 1995.
\bibitem{pop} I.~ Popovici, D.C.~ Radulescu, Characterizing the 
conformality in a Minkowski space, Ann. Inst. H. Poincar\' e 
Phys. Th\' eor. 35 (1981) 131--148. 

\bibitem{p1} P.~ \v{S}emrl,  On Hua's fundamental theorem of
 the geometry of rectangular matrices, J. Algebra 248(2002)
  366--380.

\bibitem{p2} T.~ Petek, P.~ \v{S}emrl, Adjacency preserving maps 
 on matrices and operators, Proc. R. Soc. Edinb., 
 Sect. A, Math. 132(2002)  661--684.
 
  \bibitem{p3} H.~ Radjavi, P.~ \v{S}emrl,  A short proof of Hua's fundamental 
   theorem of the geometry of hermitian matrices,
   Expos. math. 21(2003)   83--93.
   
 \bibitem{p4} P.~\v{S}emrl,  Hua's fundamental theorem of the geometry 
    of matrices and related results, Linear
    Algebra Appl. 361(2003)  161--179.
    
 \bibitem{p5} W. L.~ Chooi, M. H.~ Lim, P.~ \v{S}emrl,  Adjacency 
 preserving maps on upper triangular matrix algebras, Linear Algebra Appl.
  367(2003)  105--130.
 
 \bibitem{p6} P.~\v{S}emrl,  Hua's fundamental theorem of the geometry of matrices,
 J. Algebra 272(2004)  801--837.
 
 \bibitem{p7} L.~ Rodman, P.~\v{S}emrl, A.R.~Sourour, 
 Continuous adjacency preserving maps on real matrices,
 Can. Math. Bull. 48(2005) 267--274.
 
 \bibitem{p8} P.~\v{S}emrl. Maps on matrix and operator algebras, Jahresber.
 Deutsch. Math.-Ver. 108(2006)  91--103.
 
 \bibitem{p9} P.~\v{S}emrl. Maps on matrix spaces, Linear Algebra Appl. 413
  (2006)  364--393.

 \bibitem{w1} Z.-X. Wan, Geometry of symmetric matrices and its applications I, II, Agebra Colloq. 1
 		(1994), 97--120; 1 (1994) 201--224.
  
  \bibitem{w2} Z.-X. Wan, Geometry of matrices, Adv. stud. Pure Math. 
  24(1996) 443--453.
 
 \bibitem{w3} W.-L. Huang, On the fundamental theorems of 
 the geometries of symmetric 
 matrices,	Geom. Dedicata 78 (1999) 315--325.
 
 \bibitem{wh2} W.-L. Huang,  Adjacency preserving mappings between point-line
 geometries, Innovations in Incidence Geometry 3(2006), 25--32.
 
 \bibitem{w4} W.-L. Huang, A. E. Schroth, Adjacency preserving mappings, 
 Adv. Geom. 3(2003) 53--59.
 
\bibitem{w5} R. An, J.  Hou, L.  Zhao,
Adjacency preserving maps on the space of symmetric operators, 
Linear Algebra Appl. 405(2005) 311--324.


\end{thebibliography}
\end{document}